\documentclass[12pt]{amsart}
\usepackage{amssymb,amsfonts, amscd, amsmath}
\newtheorem{prop}{Proposition}

\newtheorem{thm}[prop]{Theorem}

\newcommand{\G}{{\Bbb G}}
\newcommand{\N}{{\Bbb N}}
\newcommand{\Q}{{\Bbb Q}}

\renewcommand{\P}{{\Bbb P}} 
\newcommand{\Z}{{\Bbb Z}}

\newcommand{\ep}{\hfill \qedsymbol \newline }

\begin{document}

\title[Dimension relation on elliptic curves]
{Relation between the dimensions of the ring generated by a vector bundle of degree zero
on an elliptic curve and a torsor trivializing this bundle}

\author{Silke Lekaus}
\address{FB 6 - Mathematik, 
Universit\"at Essen,
45117 Essen,
Germany}

\email{silke.lekaus@uni-essen.de}

\maketitle

\section{Introduction and Notations}

Let $X$ be a complete, connected, reduced scheme over a perfect field $k$.
We define Vect$(X)$ to be the set of isomorphism classes $[V]$ of vector 
bundles $V$ on $X$. We can define an addition and a multiplication on 
Vect$(X)$:
\begin{eqnarray*}
& [V]+[V']=[V\oplus V']\\
& [V] \cdot [V']=[V\otimes V'].
\end{eqnarray*}

The (naive) Grothendieck ring  $K(X)$ (see \cite{no}) is the ring associated 
to the additive monoid Vect$(X)$, that means
\[K(X)=\frac{{\Z}[\mbox{Vect}(X)]}{H},\]
where $H$ is the subgroup of ${\Z}[\mbox{Vect}(X)]$ generated by all 
elements of the form $[V\oplus V'] - [V] - [V']$.\newline
The indecomposable vector bundles on $X$ form a free basis of $K(X)$.
Since H$^0(X,\mbox{End}(V))$ is finite dimensional, the Krull-Schmidt theorem (\cite{at2}) 
holds on $X$. This means that a decomposition of a vector bundle in 
indecomposable components exists and is unique up to isomophism.\\
\indent We want to generalize a theorem of M. Nori on finite vector bundles.
A vector bundle $V$ on $X$ is called finite, if the collection $S(V)$ of all 
indecomposable components of $V^{\otimes n}$ for all integers $n\in {\Z}$ is 
finite. \\
In the following, we denote by R(V) the $\Q$-subalgebra of 
$K(X)\otimes_{\Z}{\Q}$ generated by the set $S(V)$. Thus
a vector bundle $V$ is finite if and only if
the ring $R(V)$ is of Krull dimension zero.

In \cite{no}, Nori proves the following theorem:\newline
For every finite vector bundle $V$ on $X$ there exists a finite group 
scheme $G$ and a principal $G$-bundle $\pi : P\to X$, such that $\pi^*V$ is 
trivial on $P$.
In particular, the equality
\[\dim R(V)=\dim G\, (=0)\]
holds.\\
As every vector bundle $V$ on $X$ of rank $r$ trivializes on its associated
principal GL($r$)-bundle, we can look for a group scheme $G$ of smallest 
dimension and a principal $G$-bundle on which the pullback of the vector
bundle $V$ is trivial. We might also compare the dimension of the group scheme
to dim $R(V)$.\\
In this article we consider the family of vector bundles of degree zero on 
an elliptic curve. We will prove in propositions 2 and 3 that they 
trivialize on a principal $G$-bundle with $G$ a group scheme of smallest
dimension one.\newline
As in the situation of Nori's theorem, this dimension turns out to be equal 
to the dimension of the ring $R(V)$.\\
\indent I am grateful to H\'el\`ene Esnault for suggesting the problem treated here
and for many useful discussions.

\section{Dimension relation for vector bundles of degree zero on an elliptic curve}

Let $X$ be an elliptic curve over an algebraically closed field $k$ of characteristic
zero. 
We consider vector bundles of degree zero on $X$ which
can be classified according to Atiyah (see \cite{at}).
By ${\mathcal E}(r,0)$ we denote the set of indecomposable vector bundles of 
rank $r$ and degree zero.

\begin{thm}(Atiyah \cite{at})
\begin{enumerate}
\item There exists a vector bundle $F_r\in {\mathcal E}(r,0)$, unique up to 
isomorphism, with $\Gamma (X,F_r)\neq 0$.\newline
Moreover we have an exact sequence
\begin{equation*}
\begin{array}{ccccccccc}
0 & \to & {\mathcal O}_X & \to & F_r & \to & F_{r-1} & \to & 0.
\end{array}
\end{equation*}
\item Let $E\in {\mathcal E}(r,0)$, then $E\cong L\otimes F_r$
where $L$ is a line bundle of degree zero, unique up to isomorphism (and 
such that $L^r\cong \det E$.)

\end{enumerate} 
\end{thm}

\begin{prop}\ \ 
\begin{enumerate}
\item[i)] The $\Q$-subalgebra $R(F_r)$ of $K(X)\otimes_{\Z} {\Q}$ 
generated by $S(F_r)$ is $\Q[x]$, where $x=[F_2]$, if $r$ is even, and 
$x=[F_3]$, if  $r$  is odd.
In particular, $R(F_r)$ is of Krull dimension zero.
\item[ii)]  There 
exists a principal ${\G}_a$-bundle $\pi : P \to X$ such that $\pi^*(F_r)$ is 
trivial for all $r\ge 2$.\end{enumerate}

Remark: As in Nori's case we have a correspondence of dimensions
\[\mbox{dim }R(F_r)=\mbox{dim }{\G}_a = 1.\]
\end{prop}

\noindent Proof:\newline
As proved by Atiyah in  \cite{at}, the vector bundles $F_r$ are self-dual and 
fulfill the formula
\[F_r\otimes F_s= F_{r-s+1}\oplus F_{r-s+3} \oplus \cdots \oplus 
F_{(r-s)+(2s-1)}\, \mbox{ for } s\le r.\]
For even $r$, it follows by induction that there exist integers $a_i(n)$ 
such that 
\[F_r^{\otimes n} = a_2(n)F_2\oplus a_4(n) F_4 \oplus \cdots \oplus
a_{(r-1)n -1}(n)F_{(r-1)n -1} \oplus F_{(r-1)n+1}\]
for odd $n\ge 3$, and
\[F_r^{\otimes n} =  a_1(n){\mathcal O}_X \oplus a_3(n) F_3 \oplus \cdots 
\oplus a_{(r-1)n-1}(n)F_{(r-1)n-1} \oplus F_{(r-1)n+1}\]
for even $n\ge 2$ .\\
Therefore we obtain
\[S(F_r)=\{F_i \, |\, i=1,2,3,\dots\}\,  , \mbox{ if } r \mbox{ even },\]
and $S(F_r)$ generates the subring ${\Q}[F_2]$ of $K(X)\otimes {\Q}$, because 
inductively we can write every vector bundle $F_i$ as  $p(F_2)$ for some 
polynomial $p \in \Z[x]$.

For odd $r$, Atiyah's multiplication formula gives
\[F_r^{\otimes n} =  a_1(n){\mathcal O}_X \oplus a_3(n) F_3 \oplus \cdots \oplus
a_{(r-1)n -1}(n)F_{(r-1)n -1} \oplus F_{(r-1)n+1} \]
for all $n \ge 2$.
It follows that 
\[S(F_r)=\{F_i \, |\, i \mbox{ odd }\}\, , \mbox{ if } r \mbox{ odd }.\]
For odd $r$, the set $S(F_r)$ generates the ring $R(F_r)={\Q}[F_3]$, as for odd $i$ 
each $F_i$ is $p(F_3)$ for a polynomial $p \in \Z[x]$. 

The vector bundle $F_2 $ is an element of 
$H^1(X,\mbox{GL}(2,\mathcal O))$. Because of the exact sequence
\[\begin{array}{ccccccccc}
0 & \to & {\mathcal O}_X & \to & F_2 & \to & {\mathcal O}_X & \to & 0,
\end{array}
\]
$F_2$ is even an element of $H^1(X,{\G}_a)$. Here we embed $\G_a$ into
GL($2,{\mathcal O}$) via $u \to \left(
\begin{array}{cc}
1 & u \\
0 & 1 \\
\end{array}
\right)
$.
Hence $F_2$ trivializes on a principal ${\G}_a$-bundle.
As $F_r= S^{r-1}F_2\, , r\ge 3,$ each $F_r$ trivializes on the same 
principal ${\G}_a$-bundle as $F_2$.\\
\indent As the classes $[F_r]$ are not torsion elements in $H^1(X,\mbox{GL}(2,\mathcal O))$,
none of the bundles $F_r$ can trivialize on a principal $G$-bundle with $G$ a finite 
group scheme.\ep

{\bf Remark:}
In the given examples of vector bundles $E$ there was so far not only a 
correspondence of the dimensions of the group scheme and the ring $R(E)$.
The algebra $R(E)$ was also the Hopf algebra corresponding to the group 
scheme. The following proposition shows that this is not true in general.

\begin{prop}
Let $E\cong L\otimes F_r \in {\mathcal E}(r,0)$ (see theorem 1). 
\begin{enumerate}
\item
If $L$ is not torsion, the ring $R(E)$ is isomorphic to $\Q [x,x^{-1}]
\otimes \Q [y]$ and $E$ trivializes on a principal $\G_m \times \G_a$-bundle.
\item
If $L$ is torsion, let $n\in \N$, $n\ge 1$, be the minimal number such that $L^{\otimes 
n}\cong {\mathcal O}_X$.
If $n$ and $r$ are both even, the ring $R(E)$ is isomorphic to 
\[\Q[x]/<x^{n/2} -1> \otimes\Q[y]\] 
and $E$ trivializes on a 
principal $\mu_n \times\G_a$-bundle.
There is no principal $\mu_{n/2}\times\G_a$-bundle
where $E$ is trivial.\\
If $n$ and $r$ are not both even, the ring $R(E)$ is isomorphic to
\[\Q[x]/<x^n -1> \otimes\Q[y]\]
and $E$ trivializes on a principal 
 $\mu_n\times\G_a$-bundle.
\end{enumerate}
\end{prop}

\noindent Proof:
Let $E\in {\mathcal E}(r,0)$ with $\Gamma(X,E)=0$. (If $\Gamma(X,E)\neq 0$, then
$E\cong F_r$. This case was already dealt with in proposition 2.)\\
\indent First we consider the case that $L$ is not torsion.\newline
We must distinguish between odd and even $r$. \newline
For odd $r$, Atiyah's multiplication formula ( see proof of proposition 
4) gives the following result:\\
For $m\in \N$, $m\ge 2$, the tensor power
$E^{\otimes m} \cong L^{\otimes m}\otimes F_r^{\otimes m}$
has the indecomposable components 
$ L^{\otimes m}\otimes {\mathcal O}_X, \, L^{\otimes m}\otimes 
F_3, \dots , L^{\otimes m}\otimes F_{(r-1)m +1}$,\\
the tensor power
$E^{\otimes -m} \cong L^{\otimes -m}\otimes F_r^{\otimes m}$
has the indecomposable components
$L^{\otimes -m}\otimes {\mathcal O}_X, \, L^{\otimes -m}\otimes 
F_3, \dots , L^{\otimes -m}\otimes F_{(r-1)m+1}$.\\
Thus we obtain
\[  S(E)= \left\{
\begin{array}{l}
{\mathcal O}_X, L\otimes F_r, L^{-1}\otimes F_r,\\
L^{\otimes \pm i}\otimes F_3, L^{\otimes \pm i}\otimes F_5, \dots , 
L^{\otimes \pm i}\otimes F_{(r-1)i +1}, \mbox{ i }\in \N\\
\end{array}
\right\}.
\]

The algebra $R(E)$ which is generated by $S(E)$ is the subalgebra 
of $K(X)\otimes_{\Z}\Q$ generated by $L$, $L^{-1}$ and$F_3$, thus
\[R(E)=\Q [L,L^{-1}] \otimes_{\Z} \Q[F_3].\]

For even $r$, a similar computation gives that
\[ S(E)= \left\{
\begin{array}{l}
{\mathcal O}_X, L\otimes F_r, L^{-1}\otimes F_r,\\
L^{\otimes \pm 2i}, L^{\otimes \pm 2i}\otimes F_3, \dots , L^{\otimes \pm 2i}\otimes 
F_{(r-1)2i +1},\mbox{ i }\in \N\\
L^{\otimes \pm (2i+1)}\otimes F_2, L^{\otimes \pm (2i+1)}\otimes F_4,\dots,\\ 
\mbox{  }\; \; L^{\otimes \pm (2i+1)}\otimes F_{(r-1)(2i+1) +1},\mbox{ i }\in \N\\ 
\end{array}
\right\}.\]
The ring $R(E)$, generated by $S(E)$, is the subring of 
$K(X)\otimes_{\Z} {\Q}$ which is generated by the elements $L^{\otimes 2}$, 
$L^{\otimes -2}$, $L^{-1}\otimes F_2$, therefore
\[ R(E)=\Q [L^{\otimes 2}, L^{\otimes -2}]\otimes _{\Z} \Q [L^{-1}\otimes
F_2].\] 
If $L$ is not a torsion bundle, it is clear that $L$ trivializes on a 
principal $\G_m$-bundle $P_L$. The vector bundle $E\cong L\otimes F_2$ 
trivializes on the $\G_m \times \G_a$-bundle $P_L \times_X P$, where $P$ is 
the principal $\G_a$-bundle from proposition 2,
where $F_2$ and hence all the $F_r$ trivialize.\\

Let now L be torsion and $n\in \N$, $n\ge 2$, the minimal number with 
$L^{\otimes n}\cong{\mathcal O}_X$.
As the $F_r$ are selfdual and $L^{\otimes n-1}=L^{-1}$, it suffices to 
consider positive tensor powers.\newline
Again we compute the tensor powers using Atiyah's formula to find the 
indecomposable components.

If $r$ is even and $n$ is odd, the set $S(E)$ contains the following bundles:
\[S(E)=\{ {\mathcal O}_X, L^{\otimes i}\otimes F_j \, | \, i=0,1,\dots , n-1, \,
j\in \N \}.\]
With the 
help of the multiplication formula for $F_2$ it is easy to show that all 
elements of $S(E)$ can be generated by $L$ and $F_2$. In additon, the 
relation
$L^{\otimes n}\cong {\mathcal O}_X$
holds. Hence we obtain
\[R(E)=\frac{\Q [L]}{<L^{\otimes n} -1>}\otimes_{\Z} \Q[F_2].\]

If $r$ is odd and $n$ is even or odd, the result is 
\[S(E)=\{ L^{\otimes i}\otimes F_j \, | \, i=0,1,\dots , n-1, \, j\in \N 
\mbox{ odd} \}.\]
The bundles $L$ and $F_3$ are in $S(E)$ and generate all elements of $S(E)$.
Because of the relation $L^{\otimes n}\cong {\mathcal O}_X$, the algebra $R(E)$ is
\[R(E)=\frac{\Q [L]}{<L^{\otimes n} -1>}\otimes_{\Z} \Q[F_3].\]

If $r$ and $n$ are both even
\[S(E)=\{ L^{\otimes 2i}\otimes F_{2j-1}, L^{\otimes 2i+1}\otimes F_{2j} \, 
| \, i=0,1, \dots , n/2, \, j\in \N\}.\]

The algebra R(E) is generated by $L^{\otimes 2}$ and $L\otimes F_2$. The 
generators are subject to the relation $L^{\otimes n}\cong {\mathcal O}_X$, thus
\[ R(E)= \frac{\Q[L^{\otimes 2}]}{<({L^{\otimes 2}})^{\otimes m} -1>}
\otimes \Q[L\otimes F_2],\]
where $m=n/2$.

Recall that  $n \ge 2$ is the minimal number such that 
$L^{\otimes n} \cong {\mathcal O}_X$. Thus the bundle $L$ 
trivializes on a $\mu_n$-bundle $P_L$ and not on a $\mu_m$-torsor for $m< n$.\newline
The bundle $E\cong L\otimes F_r$ then trivializes on the $\mu_n\times \G_a$-bundle
$P_L\times_X P$, where $P$ is again the principal $\G_a$-bundle from proposition 2.
We will now show that
the bundle $E$ does not trivialize on a $\mu_{n/2}\times \G_a$-bundle:\newline
If $E\cong L\otimes F_r$ trivializes on $Q\times_X P$, where $Q$ is a $\mu_m$-torsor
and $P$ a $\G_a$-torsor, then det$(L\otimes F_r) = L$ is the identity element in 
the group Pic($Q\times_X P$). But one has Pic($Q\times_X P) =$Pic$(Q)$ by homotopy invariance.
Thus $L$ must trivialize on the $\mu_m$-torsor $Q$, which is impossible for 
$m<n$. \ep

{\bf Remark:}
The correspondence between the dimension of the ``minimal'' group scheme and the 
dimension of the ring $R(E)$ also occurs in the case of vector bundles on 
the projective line, as one easily sees.\\
Let $X$ be the complex projective line ${\P}^1$ and $E:={\mathcal O} (a)$ a 
line bundle.\\
If $a=0$ we have $S(E)=\{ \mathcal O\}$ and $R(E)=Q$.\\
We define the group scheme $G$ to be $G=\mbox{Spec } {\Q}$ and the trivializing torsor 
is simply ${\P}^1$. \\
If $a\neq 0$ we can easily compute that 
$S(E)=\{{\mathcal O}(\lambda\cdot a) | \lambda \in {\Z}\}$ and
$R(E)={\Q}[x,x^{-1}]$. We define the group scheme to be 
$G={\G}_m=\mbox{Spec } {\Q}[x,x^{-1}]$.\\
The given line bundle $E$ trivializes on a principal ${\G}_m$-bundle $P_a$, 
which depends on $a$.\\
Thus we get the correspondence of 
$\dim R(E)$ and $\dim G$ in the case of a line bundle on ${\P}^1$.
This computation can easily be generalized to the case of 
vector bundles of higher rank. We illustrate this for bundles of rank two.\\
Let now $E$ be a vector bundle of rank 2 on ${\P}^1$, $E={\mathcal O}(a)\oplus 
{\mathcal O}(b)$.\\
The case $(a,b)=(0,0)$ is trivial. 
We can see at once that
$S(E)=\{{\mathcal O}\}$ and  therefore $R(E)={\Q}$.\\
The vector bundle $E$ trivializes on the principal Spec ${\Q}$ - bundle 
${\P}^1$.\\
If $(a,b)\neq (0,0)$ the computation gives that
 $S({\mathcal O}(a)\oplus {\mathcal O}(b))\, = \, S({\mathcal O}(c))$,\\
where $c=(a,b)$ (with $(a,0)=a$ and $(0,b)=b$) and
therefore  $R(E)={\Q}[x,x^{-1}]$. $E$ trivializes on the principal 
${\G}_m$-bundle $P_c$ that belongs to ${\mathcal O}(c)$ as ${\mathcal O}(a)={\mathcal 
O}(c)^\lambda$ and ${\mathcal O}(b)={\mathcal 
O}(c)^\mu$ for appropriate integers $\lambda$ and $\mu$.\\

\end{document}